\documentstyle[twoside,amstex]{article}
\input amssym.def
\input amssym.tex
\def\bc{\begin{center}}
\def\ec{\end{center}}
\def\no{\noindent}
\def\hang{\hangindent\parindent}
\def\textindent#1{\indent\llap{[#1]\enspace}\ignorespaces}
\def\re{\par\hang\textindent}
\begin{document}
\thispagestyle{empty} \vspace*{3 true cm} \pagestyle{myheadings}
\markboth {\hfill {\sl Marjan Shebani and Huanyin Chen}\hfill}
{\hfill{\sl Elementary Matrix Reduction Over J-Stable Rings
}\hfill} \vspace*{-1.5 true cm} \bc{\large\bf Elementary Matrix
Reduction Over J-Stable Rings}\ec

\vskip6mm \bc{{\bf Marjan
Sheibani}\\[2mm]
Faculty of Mathematics, Statistics and Computer Science\\
Semnan University, Semnan, Iran\\
m.sheibani1@@gmail.com}\ec

\bc{{\bf Huanyin Chen$^*$}\\[2mm]
Department of Mathematics, Hangzhou Normal University\\
Hangzhou 310036, China\\
huanyinchen@@aliyun.com}\ec

\begin{figure}[b]
\vspace{-3mm}
\rule[-2.5truemm]{5cm}{0.1truemm}\\[2mm]
{\footnotesize $^*$ Corresponding author.}
\end{figure}
\setcounter{page}{1}

\vskip10mm
\begin{minipage}{120mm}
\no {\bf Abstract:} A commutative ring $R$ is J-stable provided
that $R/aR$ has stable range 1 for all $a\not\in J(R)$. A
commutative ring $R$ in which every finitely generated ideal is
principal is called a B$\acute{e}$zout ring. A ring $R$ is an
elementary divisor ring provided that every matrix over $R$ admits
a diagonal reduction. We prove that a J-stable ring is a
B$\acute{e}$zout ring if and only if it is an elementary divisor
ring. Further, we prove that every J-stable ring is strongly completable. Various types of J-stable rings are
provided. Many known results are
thereby generalized to much wider class of rings, e.g. [2, Theorem 8], [4,
Theorem 4.1], [7, Theorem 3.7], [8, Theorem], [9, Theorem 2.1], [13, Theorem 1] and [17, Theorem 7].
\vskip3mm \no {\bf MR(2010) Subject Classification}: 13F99, 13E15,
06F20.
\end{minipage}

\vskip15mm \bc{\bf 1. INTRODUCTION}\ec

\vskip4mm \no Throughout this paper, all rings are commutative
with an identity. A matrix $A$ (not necessarily square) over a ring
$R$ admits diagonal reduction if there exist invertible matrices
$P$ and $Q$ such that $PAQ$ is a diagonal matrix $(d_{ij})$, for
which $d_{ii}$ is a divisor of $d_{(i+1)(i+1)}$ for each $i$. A
ring $R$ is called an elementary divisor ring provided that every
matrix over $R$ admits a diagonal reduction. A ring is a
B$\acute{e}$zout ring if every finitely generated ideal is
principal. Evidently, every elementary divisor ring is a
B$\acute{e}$zout ring. It is attractive to explore the conditions
under which a B$\acute{e}$zout ring (maybe with zero divisors) is an elementary divisor ring.

A ring $R$ is a Hermite ring if every $1\times 2$ matrix over $R$
admits a diagonal reduction. As is well known, a ring $R$ is
Hermite if and only if for all $a,b\in R$ there exist $a_1,,b_1\in
R$ such that $a=a_1d,b=b_1d$ and $a_1R+b_1R=R$ ([14, Theorem 1.2.5]). Thus, every
Hermite ring is a B$\acute{e}$zout ring. After Kaplansky's work on
elementary divisor rings without nonzero zero divisors, Gillman and Henriksen proved that

\vskip4mm \hspace{-1.8em} {\bf Theorem 1.1 [7, Theorem 1.1].}\ \
{\it A ring $R$ is an elementary divisor ring if and only if}
\begin{enumerate}
\item [(1)] {\it $R$ is a Hermite ring;}
\vspace{-.5mm}
\item [(2)] {\it For all $a_1,a_2,a_3\in R$, $a_1R+a_{2}R+a_{3}R=R\Longrightarrow $ there exist $p,q\in R$ such that
$pa_1R+(pa_{2}+qa_{3})R=R$.}\end{enumerate}

\vskip2mm A ring $R$ is said to have stable range 1 provided that
$aR+bR=R$ with $a,b\in R$ implies that $a+by\in R$ is invertible
for a $y\in R$. It was first introduced so as to study
stabilization in algebraic $K$-theory. Afterwards, it was studied
to deal with the cancellation problem of modules [1]. As is well
known, a regular ring $R$ is unit-regular if and only if $R$ has
stable range 1. In [3, Theorem 3], it was proved that every unit-regular
ring is an elementary divisor ring. Following McGovern, a ring $R$ has almost stable range 1
provided that every proper homomorphic image of
$R$ has stable range $1$. Evidently, every ring having stable
range 1 has almost stable range 1. Moreover, every neat ring (including FGC-domains and h-local domains) has almost stable range 1 (cf. [6]). In [7, Theorem 3.7], it is proved that
every B$\acute{e}$zout ring having almost stable range 1 is an elementary divisor
ring.

In this article, we generalize almost stable rang 1 and introduce
J-stable rings. We say that a ring $R$ is J-stable provided that
$R/aR$ has stable range 1 for all $a\not\in J(R)$. Here, $J(R)$
denote the Jacobson radical of $R$. Clearly, every ring having
almost stable 1 is J-stable. In Section 2, we shall investigate
elementary properties of J-stable rings. Various types of J-stable rings which do not have almost stable range 1 are
provided.

An element $a\in R$ is adequate if for any $b\in R$ there exist
some $r,s\in R$ such that $(1)$ $a=rs$; $(2)$ $rR+bR=R$; $(3)$
$s'R+bR\neq R$ for each non-invertible divisor $s'$ of $s$. A
B$\acute{e}$zout ring in which every nonzero element is adequate
is called an adequate ring. Kaplansky proved that for the class of
adequate domains being a Hermite ring was equivalent to being an
elementary divisor ring. This was extended to rings with
zero-divisors by Gillman and Henriksen [2, Theorem 8]. In Section 3, we
consider certain subclasses of J-stable rings by means of adequate
property. A B$\acute{e}$zout ring in which every element not in
$J(R)$ is adequate is called J-adequate. Clearly, every adequate
ring is J-adequate. For instances, regular rings and valuation
rings. A B$\acute{e}$zout ring $R$ is $\pi$-adequate provided that
for any $a\neq 0$ there exists some $n\in {\Bbb N}$ such that
$a^n\in R$ is adequate. We shall prove that every J-adequate ring
and every $\pi$-adequate ring are J-stable rings, and then enrich
the supply of such new rings by means of generalizations of
adequate rings.

Finally, we prove, in Section 4, that every J-stable ring is a
B$\acute{e}$zout ring if and only if it is an elementary divisor ring. This gives a
nontrivial generalization of [7, Theorem 3.7]. The technique here inspires us to introduce quasi adequate
rings, and prove that every quasi adequate ring is an elementary
divisor ring. This extend [13, Theorem 1] as well. Furthermore, we generalize [9, Theorem 2.1] and prove that every J-stable
ring is strongly completeable. This also extend [8, Theorem] to much wider
class of rings (maybe with zero divisors).

\vskip15mm\bc{\bf 2. $J$-Stable Rings}\ec

\vskip4mm The purpose of this section is to investigate elementary
properties of J-stable rings.

\vskip4mm \hspace{-1.8em} {\bf Theorem 2.1.}\ \ {\it Let $R$ be a
ring. Then the following are equivalent:}
\begin{enumerate}
\item [(1)] {\it $R$ is J-stable.}
\vspace{-.5mm}
\item [(2)] {\it $a_1R+a_{2}R+a_{3}R=R$ with
$a_1\not\in J(R)~\Longrightarrow $ there exists $b\in R$ such that
$a_1R+(a_{2}+a_{3}b)R=R$.} \vspace{-.5mm}
\item [(3)] {\it For any $a\not\in J(R)$, $R/a^nR$ has stable range $1$ for some $n\in {\Bbb N}$.}
\end{enumerate}
\vspace{-.5mm}  {\it Proof.}\ \ $(1)\Rightarrow (3)$ This is
obvious.

$(3)\Rightarrow (2)$ Suppose that $aR+bR+cR=R$ with $a\not\in
J(R),b,c\in R$. Then $a^nR+bR+cR=R$. By hypothesis, $R/a^nR$ has
stable range $1$. Clearly,
$\overline{b}(R/a^nR)+\overline{c}(R/a^nR)=R/a^nR$. Then we can
find some $y\in R$ such that $\overline{b+cy}\in R/a^nR$ is
invertible. Hence, $\overline{(b+cy)d}=\overline{1}$, and then
$a^nx+(b+cy)d=1$ for some $x\in R$. Therefore $aR+(b+cy)R=R$, as
desired.

$(2)\Rightarrow (1)$ Let $a\not\in J(R)$. Given
$\overline{bc}+\overline{d}=\overline{1}$ in $R/aR$, then
$ax+bc+d=1$ for some $x\in R$. By hypothesis, there exists a $y\in
R$ such that $aR+(b+dy)R=R$. Hence,
$\overline{b+dy}\big(R/aR\big)=R/aR$, and so $\overline{b+dy}\in
R/aR$ is invertible. Therefore, $R/aR$ has stable range $1$, as
desired.\hfill$\Box$

\vskip4mm \hspace{-1.8em} {\bf Corollary 2.2.}\ \ {\it Let $R$ be
a ring. then the following are equivalent:}
\begin{enumerate}
\item [(1)] {\it $R$ is J-stable.}
\vspace{-.5mm}
\item [(2)] {\it For every three elements $a,b,c\in R$ such that $a\not\in J(R)$ and $bR+cR=R$, there exists a $y\in R$ such that
$aR+(b+cy)R=R$.}
\end{enumerate}
\vspace{-.5mm}  {\it Proof.}\ \ $(1)\Rightarrow (2)$ This is
obvious.

$(2)\Rightarrow (1)$ Let $a\not\in J(R)$. Suppose that
$\overline{b}(R/aR)+\overline{c}(R/aR)=R/aR$ with $b,c\in R$.
Write $ax+by+cz=1$ for some $x,y,x\in R$. Hence, $by+(cz+ax)=1$.
By hypothesis, $aR+(b+(cz+ax)t)R=R$ for a $t\in R$. Thus,
$aR+(b+czt)R=R$, and so $\overline{b+czt}\in R/aR$ is invertible.
This implies that $R/aR$ has stable range $1$. Therefore $R$ is
J-stable, by Theorem 2.1.\hfill$\Box$

\vskip4mm \hspace{-1.8em} {\bf Corollary 2.3.}\ \ {\it Let $e$ be
an idempotent of a J-stable ring $R$, then $eRe$ is J-stable.}
\vskip2mm\hspace{-1.8em} {\it Proof.}\ \ For every three elements
$a,b,c\in eRe$ such that $a\not\in J(eRe)$ and
$b(eRe)+c(eRe)=eRe$, we see that $a\not\in J(R)$ and
$(b+1-e)R+cR=R$. Since $R$ is $J$-stable, it follows by Corollary
2.2 that $aR+(b+1-e+cy)R=R$ for a $y\in R$. Write
$ax+(b+1-e+cy)z=1$ for some $x,z\in R$. Then $(1-e)ze=0$, and so
$ze=eze$. This implies that $a(exe)+(b+c(eye))(eze)=e$. Hence,
$a(eRe)+(b+c(eye))(eRe)=eRe$. By using Corollary 2.2 again, $eRe$
is $J$-stable.\hfill$\Box$

\vskip4mm Following Moore and Steger, a ring $R$ is called a B-ring provided
that $a_1R+\cdots +a_{n}R=R (n\geq 3)$ with $(a_1,\cdots
,a_{n-2})\nsubseteq J(R)~\Longrightarrow $ there exists $b\in R$
such that $a_1R+\cdots +a_{n-2}R+(a_{n-1}+a_{n}b)R=R$. Elementary
properties of such rings have been studied in [9].
Surprisingly, we shall prove the classes of J-stable rings and
B-rings coincide with each other. That is,

\vskip4mm \hspace{-1.8em} {\bf Proposition 2.4.}\ \ {\it Let $R$ be
a ring. then the following are equivalent:}
\begin{enumerate}
\item [(1)] {\it $R$ is J-stable.}
\vspace{-.5mm}
\item [(2)] {\it $a_1R+\cdots +a_{n}R=R (n\geq 3)$ with $(a_1,\cdots
,a_{n-2})\nsubseteq J(R)~\Longrightarrow $ there exists $b\in R$
such that $a_1R+\cdots +a_{n-2}R+(a_{n-1}+a_{n}b)R=R$.}
\end{enumerate}
\vspace{-.5mm}  {\it Proof.}\ \ $(1)\Rightarrow (2)$ The assertion is true for $n=3$, by Theorem 2.1. Suppose the result holds for $n=k (k\geq 3)$. Given
$a_1R+\cdots +a_{k+1}R=R (n\geq 3)$ with $(a_1,\cdots
,a_{k-1})\nsubseteq J(R)$, then there are $x_1,\cdots ,x_{k+1}\in R$ such that $(a_1x_1+a_2x_2)+a_3x_3+\cdots +a_{k-1}x_{k-1}+a_kx_k+x_{k+1}x_{k+1}=1$. Hence,
$(a_1x_1+a_2x_2)R+a_3R+\cdots +a_{k-1}R+a_kR+a_{k+1}R=R$. By hypothesis,
$(a_1x_1+a_2x_2)R+a_3R+\cdots +a_{k-1}R+(a_k+a_{k+1}z)R=R$ for some $z\in R$. Therefore
$a_1R+a_2R+\cdots +a_{k-1}R+(a_k+a_{k+1}z)R=R$. By induction, we complete the proof.

$(2)\Rightarrow (1)$ This is obvious.\hfill$\Box$

\vskip4mm By virtue of Proposition 2.4, we see that
J-stable rings and B-rings coincide, but we prefer to use this
new concept as the preceding Theorem 2.1 shows that it is close to
stable range 1. This observation provides many class of such
rings. For instances, semi-local rings, local rings, $\pi$-regular
rings, regular rings, Noetherian rings in which every proper prime
ideal is maximal (in particular, Dedekind domain) are all
J-stable. Let $a\in R$ and $mspec(a)=\{ M\in Max(R)~|~a\in M\}$.
Further, we have

\vskip4mm \hspace{-1.8em} {\bf Example 2.5.}\ \ Let $R$ be a ring
in which $mspec(a)$ is finite for all $a\not\in J(R)$. Then $R$ is
J-stable.

\vskip4mm Let $R$ be ring in which $R/aR$ is semilocal for all
$a\not\in J(R)$. Then $R$ is J-stable. As every semilocal ring has
stable range 1, we are done by Theorem 2.1.

\vskip4mm \hspace{-1.8em} {\bf Lemma 2.6.}\ \ {\it A ring $R$ has
stable range $1$ if and only if $a_1R+a_{2}R+a_{3}R=R$ implies
that $a_1R+(a_{2}+a_{3}b)R=R$ for some $b\in
R$.}\vskip2mm\hspace{-1.8em} {\it Proof.}\ \ $\Longrightarrow$
Suppose $a_1R+a_{2}R+a_{3}R=R$. Then
$a_1x_1++a_{2}x_{2}+a_{3}x_{3}=1$; hence,
$a_{2}R+(a_1x_1+a_{3}x_{3})R=R$. Thus, we have a $y\in R$ such
that $a_{2}+(a_1x_1+a_{3}x_{3})y=u\in U(R)$. Hence,
$a_1x_1yu^{-1}+(a_{2}+a_{3}x_{3}y)u^{-1}=1$. Therefore,
$a_1R+(a_{2}+a_{3}b)R=R$, where $b=x_{3}y$.

$\Longleftarrow$ This is obvious.\hfill$\Box$

\vskip4mm \hspace{-1.8em} {\bf Theorem 2.7.}\ \ {\it Let
$\{R_{i}~|~i\in I\}~(|I|\geq 2)$ be a family of rings. Then the
direct product $R=\prod R_{i}$ of rings $R_i$ is J-stable if and
only if each $R_{i}$ has stable range $1$.}
\vskip2mm\hspace{-1.8em} {\it Proof.}\ \ $\Longrightarrow$ Given
$aR_1+bR_1+cR_1=R_1$ with $a,b,c\in R$, then $$(a,1,1,\cdots )R+(b,0,0,\cdots )R+(c,0,0,\cdots )R=R.$$ Clearly,
$(a,1,1,\cdots, 1)\not\in J(R)$. By hypothesis, there exists
$(y,y_2,y_3,\cdots )\in R$ such that $$(a,1,1,\cdots )R+\big((b,0,0,\cdots )+(c,0,0,\cdots )(y,y_2,y_3,\cdots )\big)R=R.$$ Therefore, $aR_1+(b+cy)R_1=R_1$.
In light of Lemma 2.6, $R_1$ has stable range $1$. Likewise, $R_i$
has stable range $1$ for $i\neq 1$. Therefore each $R_{i}$ has
stable range $1$.

$\Longleftarrow $ Since each $R_i$ has stable range one, we see that $R=\prod R_{i}$ has stable range $1$.
Therefore, $R$ is J-stable.\hfill$\Box$

\vskip4mm \hspace{-1.8em} {\bf Corollary 2.8.}\ \ Let
$L=\prod\limits_{i\in I}R_{i}$ be the direct product of rings
$R_i\cong R$ and $|I|\geq 2$. Then $L$ is J-stable if and only if
$R$ has stable range $1$ if and only if $L$ has stable range $1$.

\vskip4mm Thus, we see that ${\Bbb Z}\times {\Bbb Z}$ is not
J-stable, while ${\Bbb Z}$ is J-stable.

\vskip4mm \hspace{-1.8em} {\bf Proposition 2.9.}\ \ {\it Let $R$ be a
ring, and let $I$ be an ideal of $R$. If $I\subseteq J(R)$, then
the following are equivalent:}
\begin{enumerate}
\item [(1)] {\it $R$ is J-stable.}
\vspace{-.5mm}
\item [(2)] {\it $R/I$ is J-stable.}
\end{enumerate}
\vspace{-.5mm}  {\it Proof.}\ \  $(1)\Rightarrow (2)$ Let $R$ be a
J-stable ring and let $\bar {a_1}\overline{R}+\bar
{a_2}\overline{R}+\bar{a_3}\overline{R}= \overline{R}$, where
$\bar{a_1}\not\in J(\overline{R})$,  so there are $\bar{r_1},
\bar{r_2}, \bar{r_3}\in\overline{R}$ such that
$\bar{a_1}\bar{r_1}+\bar{a_2}\bar{r_2}+\bar{a_3}\bar{r_3}
=\bar{1}$. Hence, $a_1r_1+a_2r_2+a_3r_3=1+x$ for some $x\in I$ as
$I\subseteq J(R)$, $1+x$ is a unit and then
$(a_1r_1+a_2r_2+a_3r_3)R=R$, hence $a_1R+a_2R+a_3R=R$ with
$a_1\not\in J(R)$. By Theorem 2.1, there exists $b\in R$ such that
$a_1R+(a_2+a_3b)R=R$ and then
$\bar{a_1}\overline{R}+(\overline{a_2+a_3b})
\overline{R}=\overline{R}$.

$(2)\Rightarrow (1)$ Let $a_1R+a_2R+a_3R=R$, with $a_1\not\in
J(R)$, so we have $\bar {a_1}\overline{R}+\bar
{a_2}\overline{R}+\bar{a_3}\overline{R}=\overline{R}=R/I$ as $R/I$
is a J-stable ring, there exists $\bar{b}\in\overline{R}$ such
that $\bar{a_1}\overline{R}+(\overline{a_2+a_3b})
\overline{R}=\overline{R}$ so $a_1r_1+(a_2+a_3b)r_2=1+y$ for some
$y\in I$. As $1+y\in U(R)$, we have $a_1R+(a_2+a_3b)R=R$. According to Theorem 2.1, $R$ is J-stable.\hfill$\Box$

\vskip4mm It follows by Proposition 2.9 that $R$ is J-stable if and
only if $R/J(R)$ is J-stable. Also we see that every homomorphic
image of a J-stable ring is J-stable. Thus, $R$ is J-stable if and
only if $R/aR$ is J-stable for any $a\in R$.

\vskip4mm \hspace{-1.8em} {\bf Corollary 2.10.}\ \ {\it Let $R$ be
a ring. Then the following are equivalent:}
\begin{enumerate}
\item [(1)] {\it $R$ is J-stable.}
\vspace{-.5mm}
\item [(2)] {\it $R[[x_1,\cdots ,x_n]]$ is J-stable.}
\end{enumerate}
\vspace{-.5mm}  {\it Proof.} Let $\psi: R[[x_1,\cdots
,x_n]]\longrightarrow R$ is defined by $\psi(f(x_1,\cdots
,x_n))=f(0,\cdots ,0)$ it is easily prove that $\psi$ is a ring
epimorphism such that $ker\psi\subseteq J(R[[x_1,\cdots ,x]])$,
now the result follows from the Proposition 2.9.\hfill$\Box$

\vskip4mm Let $R$ be a ring, and let $M$ be an $R$-$R$-bimodule.
Then the trivial extension $T(R,M)$ is the ring $\{ (r,m)~|~r\in
R, m\in M\}$, where the operations are defined as follows: For any
$r_1,r_2\in R, m_1, m_2\in M$,
$$\begin{array}{c}
(r_1,m_1)+(r_2,m_2)=(r_1+r_2, m_2+m_2),\\
(r_1,m_1)(r_2,m_2)=(r_1r_2, r_1m_2 + m_1r_2).
\end{array}$$

\vskip2mm \hspace{-1.8em} {\bf Corollary 2.11.}\ \ {\it Let $R$ be
a ring, and let $M$ be an $R$-$R$-bimodule. Then the following are
equivalent:}
\begin{enumerate}
\item [(1)] {\it $R$ is J-stable.}
\vspace{-.5mm}
\item [(2)] {\it $T(R,M)$ is J-stable.}
\end{enumerate}
\vspace{-.5mm}  {\it Proof.}\ \ Let $\psi:R\longrightarrow
T(R,M)/J(T(R,M)$ be such that $\psi(r)=(\overline{r,0)}$ for any
$r\in R$ $Ker(\psi)=\lbrace r\in R \vert  (r,0)\in
J(T(R,M))\rbrace=\lbrace r\in R \vert r\in J(R)\rbrace$. As
$J(T(R,M))=\lbrace (r,m) \vert r\in J(R), m\in M\rbrace$. Also for
any $\overline{(r,m)} \in T(R,M)/J(T(R,M)$ we can write
$\overline{(r,m)}=\overline{(r,0)}+\overline{(0,m)}=\psi
(r)+0=\psi(r)$ that shows $\psi$ is surjective, now we can get the
result by Proposition 2.9.\hfill$\Box$

\vskip4mm As an immediate consequence, we deduce that a ring $R$
is J-stable if and only if the ring $\{ \left(\begin{array}{cc}
a&b\\
0&a
\end{array}
\right)~|~a,b\in R\}$ is J-stable.

\vskip4mm Clearly, every ring of stable range 1 is J-stable.
Further, all rings having almost stable rang 1 is J-stable. But
the J-stable rings having not almost stable range 1 are rich.

\vskip4mm \hspace{-1.8em} {\bf Example 2.12.}\ \ {\it Let $R=\{
a+bx~|~a,b\in {\Bbb Z}, x^2=0\}$. Then $R$ is a J-stable ring,
while it does not have almost stable range
1.}\vskip2mm\hspace{-1.8em} {\it Proof.}\ \ Clearly, $J(R)=\{
bx~|~a,b\in {\Bbb Z}, x^2=0\}$. Let $z\not\in J(R)$. Write
$z=c+dx, c,d\in {\Bbb Z}$ and $c\neq 0$. Construct a map $\varphi:
R/zR\to {\Bbb Z}/c{\Bbb Z}, \overline{a+bx}\mapsto \overline{a}.$
Then $(R/zR)/ker(\varphi)\cong {\Bbb Z}/c{\Bbb Z}$. In view of
Example 2.5, ${\Bbb Z}/c{\Bbb Z}$ has stable range 1. On the other
hand, $ker(\varphi)=\{ \overline{bx}~|~b\in R\}\subseteq
J\big(R/zR\big)$. Thus, $R/zR$ has stable range 1. This
implies that $R$ is J-stable. Obviously, $R/xR\cong {\Bbb Z}$ does not have stable range 1.
Therefore $R$ does not have almost stable range 1.\hfill$\Box$

\vskip4mm \hspace{-1.8em} {\bf Example 2.13.}\ \ Let $R= \{ z_0 +
a_1x+a_2x^2 + \cdots +~|~z_0\in {\Bbb Z}, a_i\in {\Bbb Q}\}$. Then
$R$ is a J-stable ring, while it does not have almost stable range
1. Here, $J(R)=x{\Bbb Q}[[x]]$, and then $R/J(R)\cong {\Bbb Z}$
does not have stable range 1. Thus, $R$ does not have almost
stable range 1. Let $f(x)=z_0+a_1x+a_2x^2+\cdots \not\in J(R)$.
Then $z_0\neq 0$. Let $I=\big(f(x)\big)$. Then
$f(x)=z_0\big(1+z_0^{-1}a_1x+z_0^{-1}a_2x^2+\cdots \big)\in I$.
This implies that $z_0\in I$. For any $g(x)\in J(R)$, we see that
$g(x)\in x{\Bbb Q}[[x]]$. Clearly, $z_0\big(1+z_0^{-1}g(x)\big)\in
I$, and so $z_0+g(x)\in I$. This implies that $g(x)\in I$. Thus,
$J(R)\subseteq I$. We infer that $I/J(R)$ is a proper ideal of
$R/J(R)$. Since $R/I\cong R/J(R)/I/J(R)$, we see that $R/I$ has
stable range 1, therefore $R$ is J-stable.

\vskip4mm \hspace{-1.8em} {\bf Example 2.14.}\ \ Let $R$ be the
collection of all elements of the form $a\alpha+b\beta+c\gamma+d$,
with $a,b,c,d\in {\Bbb Z}_2$, where $\alpha,\beta,\gamma$ satisfy
the relations
$$\alpha^2=\beta^2=\gamma^2=\alpha\beta=\beta\alpha=\alpha\gamma=\gamma\alpha=\beta\gamma=\gamma\beta=0.$$ Then $R[x]$ is J-stable, but $R[x]$ does not have almost stable range 1.
As $R$ is consisted entirely of nilpotent elements $\{
0,\alpha,\beta,\gamma,\alpha+\beta,\alpha+\beta,\beta+\gamma,\alpha+\beta+\gamma$
and units
$1,1+\alpha,1+\beta,1+\gamma,1+\alpha+\beta,1+\alpha+\gamma,1+\beta+\gamma,1+\alpha+\beta+\gamma$.
It follows by [11, Example 3.4] and Proposition 2.4 that $R[x]$ is J-stable. Let
$I=(\alpha)[x]$. Then $I$ is a proper ideal of $R[x]$. Since
$R[x]/I\cong \big(R/(\alpha)\big)[x]$ has not stable range 1.
Therefore, $R[x]$ does not have almost stable range 1.

\vskip4mm By [9, Theorem 2.7] and Proposition 2.4, we see that
$R[x,y]=R[x][y]$ is not J-stable for an arbitrary ring $R$.

\vskip4mm \hspace{-1.8em} {\bf Example 2.15.}\ \ Let $K$ be a
field and $x$ be an undeterminate on $K$, and let $R=K[x]$. As $R$ is a principal ideal domain, it has almost stable range 1,
which does not have stable range 1. So it is J-stable by Theorem 2.1. It follows from Corollary 2.10 that
$R[[y]]=K[x][[y]]$ is J-stable. Now we have $J(R[[y]])=(y, J(R))$
which is non-zero ideal of $R[[y]$. Clearly, $R[[y]]/J(R[[y]])
\cong R/J(R)$. As $R/J(R)$ does not have stable range 1, and so
$R[[y]]/J(R[[y]])$ does not have stable rang 1. Therefore it shows that
$R[[y]]$ does not have almost stable range 1, and then we are
through.

\vskip15mm\bc{\bf 3. Certain Subclasses}\ec

\vskip4mm An element $a\in R$ is clean if it is the sum of an idempotent and a unit. A ring $R$ is clean provided that every element in $R$
is clean. By virtue of [1, Theorem 17.2.2], every
clean ring has stable range 1.

\vskip4mm \hspace{-1.8em} {\bf Lemma 3.1 [15, Theorem 2 and Theorem 4].}\ \ {\it Let $R$
be a B$\acute{e}$zout ring. If $a\in R$ is adequate, then $R/aR$
is clean.}

\vskip4mm \hspace{-1.8em} {\bf Theorem 3.2.}\ \ {\it Every
J-adequate ring is J-stable.} \vskip2mm \hspace{-1.5em}{\it
Proof.}\ \ Let $R$ be J-adequate, and let $a\not\in J(R)$. By
hypothesis, $a\in R$ is adequate. In view of Lemma 3.1, $R/aR$ is
clean, and so $R/aR$ has stable range 1. Therefore $R$ is J-stable,
as asserted.\hfill$\Box$

\vskip4mm Recall that a ring $R$ is called an NJ-ring if every
element $a\not\in J(R)$ is regular [10]. For instance, every
regular ring and every local ring are NJ-rings.

\vskip4mm \hspace{-1.8em} {\bf Example 3.3.}\ \ {\it Every
B$\acute{e}$zout NJ-ring is J-adequate.} \vskip2mm\hspace{-1.8em}
{\it Proof.}\ \ Let $R$ be a B$\acute{e}$zout NJ-ring. In view of
[10, Theorem 2], R must be a regular ring, a local ring or
isomorphic to the ring of a Morita context with zero pairings
where the underlying rings are both division ring. If R is
regular, then it is adequate. If $R$ is local, it is J-adequate.
If $R$ is isomorphic to the ring of a Morita context $T =
(A,B,M,N,\varphi,\phi)$ with zero pairings, where the underlying
rings are division rings $A$ and $B$. Then $R$ is not commutative,
a contradiction. Therefore $R$ is J-adequate.\hfill$\Box$

\vskip4mm It should be noted that adequate rings were studied in many papers (see [15]). Obviously, every adequate ring is J-adequate. But the converse is not true, as the following shows.

\vskip4mm \hspace{-1.8em} {\bf Example 3.4.}\ \ {\it Let $R= \{ z_0 +
a_1x+a_2x^2 + \cdots +~|~z_0\in {\Bbb Z}, a_i\in {\Bbb Q}\}$. Then
$R$ is a J-adequate ring, while it is not an adequate ring.} \vskip2mm\hspace{-1.8em}
{\it Proof.}\ \ As in [12, Example 3.3], $R$ is a B$\acute{e}$zout domain, but it is not an adequate ring.
Let $f(x)=y+b_1x+b_2x^2+\cdots\not\in J(R)$ and $g(x)=z+c_1x+c_2x^2+\cdots \in R$. Then $y\neq 0$.
Since ${\Bbb Z}$ is a principal ideal domain, it is adequate. Thus, there exist $s,t\in R$ such that
$y=st, (s,z)=1$, and that $(t',z)\neq 1$ for any non-unit divisor $t'$ of $t$.
If $(s,t)\neq 1$, then we have a nonunit $d\in R$ such that $(s,t)=d$. Hence,
$(d,z)\neq 1$, and then $(s,z)\neq 1$, an absurd. Therefore $(s,t)=1$, and so
$f(x)=\big(s+d_1x+d_2x^2+\cdots )\big(t+e_1x+e_2x^2+\cdots )$, where $d_i$ and $e_i$ are solutions of the equations:
$$\begin{array}{c}
se_1+d_1t=b_1;\\
se_2+d_2t=b_2-d_1e_1;\\
se_3+d_3t=b_3-d_1e_2-d_2e_1;\\
\vdots
\end{array}$$ Set $s(x)=s+d_1x+d_2x^2+\cdots $ and $t(x)=t+e_1x+e_2x^2+\cdots $.
Clearly, we can find some $k,l\in {\Bbb Z}$ such that $ks+lz=1$. Hence, $1-\big(ks(x)+lg(x)\big)\in J(R)$. Thus,
$ks(x)+lg(x)\in U(R)$. This shows that $\big(s(x),g(x))=1$. If $t'(x)=m+f_1x+f_2x^2+\cdots $ is a nonunit divisor of $t(x)$, then $m$ is a nonunit divisor of $t$. By hypothesis, $(m,z)\neq 1$. This implies that $\big(t'(x),g(x)\big)\neq 1$. Therefore $R$ is J-adequate, as asserted.
\hfill$\Box$

\vskip4mm \hspace{-1.8em} {\bf Lemma 3.5.}\ \ {\it Let $R$ be
a J-adequate ring. Then every prime ideal not in $J(R)$ contains
in a unique maximal ideal of $R$.} \vskip2mm \hspace{-1.5em}{\it
Proof.}\ \ Let $P\nsubseteq J(R)$ be a prime ideal of $R$. Then
$P$ contains at least one adequate element.
As in the proof of [14, Proposition 3.2.4], $P$ is contained in a unique
maximal ideal of $R$.\hfill$\Box$

\vskip4mm Let $R$ be a J-adequate ring. Then $R/P$ is a local
ring for all prime ideal $P\nsubseteq J(R)$. Recall that a ring is
a pm ring provided that every prime ideal is contained in a
unique maximal ideal of $R$ (see [5]). As is well known, a ring $R$ is a pm ring if and only if $a+b=1$ implies that $(1-ax)(1-by)=0$ for some $x,y\in R$. Further, we derive

\vskip4mm \hspace{-1.8em} {\bf Theorem 3.6.}\ \ {\it Let $R$ be a
J-adequate ring. Then $R/J(R)$ is a pm ring.} \vskip2mm
\hspace{-1.5em}{\it Proof.}\ \ Let $a,b\in R$ be such that
$a+b=1$. Set $A=\{ 1+ax~|~x\in R\},B=\{ 1+bx~|~x\in R\}$ and
$T=AB$. We claim that $J(R)\bigcap T\neq \emptyset$. If not, we
have a nonempty set $\Omega=\{ I\lhd R~|~I\bigcap T=\emptyset\}$,
as $J(R)\in \Omega$. Given $I_1\subseteq I_2\subseteq \cdots $ in
$\Omega$, then $I_i\bigcap T=\emptyset$. Hence,
$\big(\bigcup\limits_{i}I_i\big)\bigcap
T=\bigcup\limits_{i}\big(I_i\bigcap T\big)=\emptyset$, and so
$\bigcup\limits_{i}I_i\in \Omega$. Thus, $\Omega$ is inductive. By
Zorn's Lemma, there exists an ideal $P$ of $R$, which is maximal
in $\Omega$. If $P\not\in Spec(R)$, there exist $c,d\in R$ such
that $c,d\not\in P$, while $cd\in P$. Then $(RcR+P)\bigcap T,
(RdR+P)\bigcap T\neq \emptyset$. It follows that
$\big((RcR+P)(RdR+P)\big)\bigcap T\neq \emptyset$, and so
$(RcdR+P)\bigcap T\neq \emptyset$, a contradiction. Therefore
$P\in Spec(R)$.

If $RaR+P=R$, then $az+p=1$ for some $z\in R,p\in P$. This implies
that $p=1-az\in A\subseteq T$, an absurd. Thus, $RaR+P\subseteq M$
for an $M\in Max(R)$. Likewise, $RbR+P\subseteq N$ for an $N\in
Max(R)$. Hence, $P\subseteq M\bigcap N$. But $M\neq N$; otherwise,
$1=a+b\in M=N$. By virtue of Lemma 3.5, $P\subseteq J(R)$, and so
$1-az\in J(R)$. It follows that $J(R)\bigcap T\neq \emptyset$, a
contradiction. Accordingly, $J(R)\bigcap (AB)\neq \emptyset$.
Therefore $(1+ar)(1+bs)\in J(R)$ for some $r,s\in R$.

Given $\overline{a+b}=\overline{1}$ in $R/J(R)$, then $ax+by=1$
for some $x,y\in R$. By the preceding discussion,
$(1+axr)(1+bys)\in J(R)$ for some $r,s\in R$. Hence,
$\overline{(1+axr)(1+bys)}=\overline{0}$. Therefore $R/J(R)$ is
a pm ring, as asserted.\hfill$\Box$

\vskip4mm \hspace{-1.8em} {\bf Corollary 3.7.}\ \ {\it Let $R$ be a J-adequate ring. Then the following are equivalent:}
\begin{enumerate}
\item [(1)] {\it $Max(R)$ is zero-dimensional.}
\vspace{-.5mm}
\item [(2)] {\it $Max(R/J(R))$ is zero-dimensional.}
\end{enumerate}
\vspace{-.5mm}  {\it Proof.}\ \ Since $R$ is J-adequate, $R/J(R)$ is a pm ring by Theorem 3.6. In light of [1, Corollary 17.1.14], $R/J(R)$ is clean if and only if $Max(R/J(R))$ is zero-dimensional. We note that $Max(R)$ is zero-dimensional if and only if $R/J(R)$ is clean, and therefore the result follows.\hfill$\Box$

\vskip4mm Our next aim is to explore $\pi$-adequate rings. This
will provide a new type of J-stable rings. We are now ready to
prove:

\vskip4mm \hspace{-1.8em} {\bf Theorem 3.8.}\ \ {\it Every
$\pi$-adequate ring is J-stable.}\vskip2mm \hspace{-1.5em}{\it
Proof.}\ \ Let $R$ be a $\pi$-adequate ring, and let $a\neq 0$.
Then $a^n\in R$ is adequate for some $n\in {\Bbb N}$. In view of
Lemma 3.1, $R/a^nR$ is clean. By virtue of [1, Theorem 17.2.2], $R/a^nR$ has
stable range 1. Therefore $R$ is J-stable, by Theorem
2.1.\hfill$\Box$

\vskip4mm We now study the possible transfer of the $\pi$-adequate
property to homomorphic images and direct products.

\vskip4mm \hspace{-1.8em} {\bf Proposition 3.9.}\ \ {\it If $R$ is
a $\pi$-adequate ring and $I$ is an ideal of $R$ contained in
$J(R)$, then $R/I$ is a $\pi$-adequate ring.}\vskip2mm
\hspace{-1.5em}{\it Proof.}\ \ Let $\overline{a}, \overline{b}\in
R/I$ with $\overline{a}\neq \overline{0}$. Then $a,b\in R$ and
$a\neq 0$. By hypothesis, there exists some $n\in {\Bbb N}$ such
that $a^n\in R$ is adequate. Thus, $a^n=rs$ and $(r,b)=1$. Hence,
$\overline{a}^n=\overline{rs}$ and
$(\overline{r},\overline{b})=\overline{1}$. Let $\overline{s'}$ be
a nonunit in $R/I$ which divides $\overline{s}$. Then $s'$ is a
nonunit in $R$. Further, $s'$ divides $s+k$ for some $k\in J(R)$.
Thus, $(s',s)\neq 1$. Since $R$ is B$\acute{e}$zout, write
$(s',s)=(u)$. Clearly, $u\in R$ is a nonunit. As $u~|~s$, we see
that $(u,b)\neq 1$. Hence, $(\overline{u},\overline{b})\neq 1$. It
follows from $\overline{u}~|~\overline{s'}$ that
$(\overline{s'},\overline{b})\neq 1$. This completes the
proof.\hfill$\Box$

\vskip4mm \hspace{-1.8em} {\bf Lemma 3.10.}\ \ {\it Let $a\in R$
be adequate. Then $a^n\in R$ is adequate for all $n\in {\Bbb
N}$.}\vskip2mm \hspace{-1.5em}{\it Proof.}\ \ This is obvious as
in the proof of [14, Proposition 3.2.2].\hfill$\Box$

\vskip4mm \hspace{-1.8em} {\bf Theorem 3.11.}\ \ {\it Let
$\{R_{i}~|~i\in I\}~(2\leq |I|<\infty)$ be a family of rings. Then
the direct product $R=\prod R_{i}$ of rings $R_i$ is
$\pi$-adequate if and only if}\begin{enumerate}
\item [(1)]{\it each $R_{i}$ is $\pi$-adequate;}
\vspace{-.5mm}
\item [(2)]{\it $0\in R_i$ is adequate for all $i\in I$.}
\end{enumerate}
\vspace{-.5mm}  {\it Proof.}\ \ $\Longrightarrow$ (1) Let $a_1,
b_1\in R_1$. Then $(a_1,1,0,\cdots,0),(b_1,1,0,\cdots ,0)\in R$.
By hypothesis, $(a_1,1,0,\cdots, 0)^n\in R$ is adequate for some
$n\in {\Bbb N}$. Hence, there exists some $(r_1,r_2,r_3,\cdots
,r_m)$, $(s_1,s_2,s_3,\cdots ,s_m)\in R$ such that
$(a_1,1,0,\cdots, 0)^n =(r_1,r_2,r_3,\cdots ,r_m)$,
$(s_1,s_2,s_3,\cdots ,s_m)$ and $\big((r_1,r_2,r_3,\cdots ,r_m),
(b_1,1,0,\cdots ,0)\big)=(1,1,1,\cdots ,1)$.
Then we have $a_1^n=r_1s_1$ and $(r_1,b_1)=1$. Now suppose that $s'_1$ be a nonunit divisor of $s_1$, then
$(s'_1,1,1,\cdots, 1)$ is a nonunit divisor of $(s_1,s_2,s_3,\cdots,s_m)$. Now assume that $(s'_1,b_1)=1$ then
$\big((s'_1,1,1,\cdots,1),(b_1,1,0,\cdots, 0)\big)=(1,1,1,\cdots ,1)$, which is a contradiction. Thus, $R_1$ is a $\pi$-adequate ring.
Likewise, each $R_i (i\geq 2)$ is $\pi$-adequate.

(2) Choose $(0,1,0,\cdots ,0)\in R$. Let $b_1\in R$. Then there
exists some $(r_1,r_2,r_3,\cdots ,r_m)$, $(s_1,s_2,s_3,\cdots
,s_m)\in R$ such that $(0,1,0,\cdots, 0)^n =(r_1,r_2,r_3,\cdots
,r_m)$, $(s_1,s_2,s_3,\cdots ,s_m)$ and $\big((r_1,r_2,r_3,\cdots
,r_m), (b_1,1,0,\cdots ,0)\big)=(1,1,1,\cdots ,1)$. Hence,
$0=r_1s_1, (r_1,b_1)=1$. If $s_1'$ is a nonunit divisor of $s_1$.
Then $(s_1',1,\cdots ,1)$ is a nonunit divisor of
$(s_1,s_2,s_3,\cdots ,s_m)$. By hypothesis, $\big((s_1',1,\cdots
,1), (b_1,1,0,\cdots ,0)\big)\neq (1,1,\cdots ,1)$. It follows
that $(s_1',b_1)\neq 1$. Therefore $0\in R_1$ is adequate.
Similarly, $e_i\in R_i (i\geq 2)$ is adequate.

$\Longleftarrow $ Let $(a_1,a_2,\cdots, a_n), (b_1,b_2,\cdots,b_n)
\in R$. By hypothesis, every element in $R_i$ is $\pi$-adequate.
Then for each $1\leq i\leq n, a_i\in R_i$ is $\pi$-adequate, so
there exist $n_i\in {\Bbb N}$ such that $a_i^{n_i}$ is an adequate
element of $R_i$. Set $m=\prod\limits_{i=1}^{n}n_i$. By virtue of
Lemma 3.10, $a_i^{m}\in R_i$ is an adequate element for each
$1\leq i\leq n$. Thus, there are $r_i, s_i\in R_i$ such that
$a_i^m=r_is_i$ and $(r_i,b_i)=1$ and for each nonunit divisor
$s'_i$ of $s_i$, $(s'_i,b_i)\neq 1$. Then $(r_1,r_2,\cdots
,r_n)(s_1,s_2,\cdots,s_n)=(a_1^m,a_2^m,\cdots,a_n^m)=(a_1,a_2,\cdots,a_n)^m$.
Now let $(s'_1,s'_2,\cdots, s'_n)$ be a nonunit divisor of
$(s_1,s_2,\cdots,s_n)$, then there exists some $s'_j$ which is a
nonunit divisor of $s_j$. Thus, $(s'_j,b_j)\neq 1$. This implies
that $\big((s'_1,s'_2,\cdots ,s'_n),(b_1,b_2,\cdots,b_n)\big)\neq
1.$ Therefore $R$ is an adequate ring.\hfill$\Box$

\vskip4mm Let $R$ be a domain. We note that $R$
is a $\pi$-adequate ring if and only if $R$ is an adequate ring.
One direction is obvious. As in the proof of [14, Proposition 3.2.3], every divisor
of an adequate element in a B$\acute{e}$zout domain is adequate,
and then the converse is true.

As is well known, every regular ring is adequate. Furthermore, we derive

\vskip4mm \hspace{-1.8em} {\bf Theorem 3.12.}\ \ {\it Every
B$\acute{e}$zout ring in which every prime ideal is
maximal is $\pi$-adequate.}\vskip2mm \hspace{-1.5em}{\it Proof.}\
\ Let $R$ be a B$\acute{e}$zout ring in which every
prime ideal is maximal. Then $R$ is $\pi$-regular. Let $a\in R$. Then we have some $m\in {\Bbb N}$
such that $a^m\in R$ is unit-regular. We claim that $a^m$ is
adequate. Let $b\in R$ be an arbitrary element so there exists
some $n\in {\Bbb N}$ such that $b^n$ is also unit regular, then
there are invertible elements, $u, v\in R$ such that $a^m=a^mua^m$
and $b^n=b^nvb^n$. Set $ e=a^mu$ and $f=b^nv$ then $e, f$ are
idempotents. Now define $r=e+f-ef$. We see that $(r)=(e,f)$, and
let $s=1-f+ef$ then $e=sr$ and $(s,f)=1$. Thus, $a^m=s(ru^{-1})$. Furthermore,
$(s,b^n)=1$, and so $(s,b)=1$. Since $r$ divides $f$,
for every non-invertible divisor $x$ of $ru^{-1}$, we get $(x,f)\neq 1$. This implies that
$(x,b^n)\neq 1$, and so $(x,b)\neq 1$. Therefore $a^m\in R$ is adequate, as required.\hfill$\Box$

\vskip4mm \hspace{-1.8em} {\bf Corollary 3.13.}\ \ {\it Every
finite B$\acute{e}$zout ring is $\pi$-adequate.}\vskip2mm \hspace{-1.5em}{\it Proof.}\
\ Since $R$ is finite, it is $\pi$-regular, and then every prime ideal of $R$ is maximal. In light of Theorem 3.12, we complete the proof.
\hfill$\Box$

\vskip10mm\bc{\bf 4. Matrices over J-Stable Rings}\ec  \vskip4mm A
ring $R$ has stable range $2$ provided that
$a_1R+a_{2}R+a_{3}R=R\Longrightarrow $ there exist $y_1,y_2\in R$
such that $\big(a_1+a_3y_1\big)R+\big(a_{2}+a_{3}y_2\big)R=R$.

\vskip4mm \hspace{-1.8em} {\bf Proposition 4.1.}\ \ {\it Every
J-stable ring has stable range 2.} \vskip2mm\hspace{-1.8em} {\it
Proof.}\ \ Let $R$ be J-stable. Suppose that $a_1R+a_2R+a_3R=R$
with $a_1,a_2,a_3\in R$. If $a_1\not\in J(R)$, then there exists
$b\in R$ such that $a_1R+(a_{2}+a_{3}b)R=R$. If $a_1\in J(R)$,
then there exist $x_1,x_2,x_3\in R$ such that
$a_2x_2+a_3x_3=1-a_1x_1\in U(R)$. Hence,
$a_2x_2(1-a_1x_1)^{-1}+a_3x_3(1-a_1x_1)^{-1}=1.$ Thus,
$$a_1+a_3x_3(1-a_1x_1)^{-1}+a_2x_2(1-a_1x_1)^{-1}=1+a_1\in U(R).$$
It follows that
$$\big(a_1+a_3x_3(1-a_1x_1)^{-1}\big)(1+a_1)^{-1}+\big(a_2+a_3\cdot
0\big)(1+a_1)^{-1}=1.$$ Hence, $\big(a_1+a_3x_3(1-a_1x_1)^{-1}\big)R+\big(a_2+a_3\cdot
0\big)R=R,$ and therefore $R$ has stable range $2$.
\hfill$\Box$

\vskip4mm But the converse of Proposition 4. 1 is not true.
For instance, ${\Bbb Z}_6[x]$ has stable range 2, but it is not J-stable, as ${\Bbb Z}_6[x]/(2)\cong {\Bbb Z}_2[x]$ has not stable range 1.
Thus, we see that $ \{ ~\mbox{rings having almost stable range 1}~\}\subsetneq \{ ~\mbox{J-stable rings}~\}\subsetneq \{ ~\mbox{rings having stable range 2}~\}.$
A ring $R$ is called completable provided that $a_1R+\cdots
+a_nR=R, a_i\in R, i=2,\cdots ,n, $ implies there is a matrix over
$R$ with first row $a_1,\cdots ,a_n$ and $det(A)=1$. As is well known, a ring $R$ is completable if and only if every
stable free $R$-module $P$, i.e. $P\oplus R^m\cong R^n$ for some $m,n\in {\Bbb N}$, is free.
Since every ring having stable range 2 is completable (cf. [11, Corollary 2.1]), we see that every
stable free module over J-stable rings is free. Moreover, we have

\vskip4mm \hspace{-1.8em} {\bf Corollary 4.2.}\ \ {\it Let $R$ be a J-stable ring. Then for any idempotent $e\in R, e\in
a_1R+\cdots +a_nR, a_i\in R, i=2,\cdots ,n, $ implies there is a
matrix over $R$ with first row $a_1,\cdots ,a_n$ and $det(A)=e$.}
\vskip2mm\hspace{-1.8em} {\it Proof.}\ \ Write $e=a_1x_1+\cdots +a_nx_n$. Then
$e=(ea_1)(ex_1)+\cdots +(ea_n)(ex_n)$. In view of Corollary 2.3, $eRe$ is a J-stable ring.
Thus, we can find a matrix $(a_{ij})\in M_n(eRe)$ whose first row is $(ea_1,\cdots ,ea_n)$ such that $det(a_{ij})=e$.
Hence, we get $\left|
\begin{array}{ccc}
a_1&\cdots&a_n\\
a_{21}&\cdots&a_{2n}\\
\vdots&\ddots&\vdots\\
a_{n1}&\cdots&a_{nn}
\end{array}
\right|=\left|
\begin{array}{ccc}
ea_1&\cdots&ea_n\\
a_{21}&\cdots&a_{2n}\\
\vdots&\ddots&\vdots\\
a_{n1}&\cdots&a_{nn}
\end{array}
\right|=e,$ as required.\hfill$\Box$

\vskip4mm We turn now to the proof of our main result.

\vskip4mm \hspace{-1.8em} {\bf Theorem 4.3.}\ \ {\it Let $R$ be a
J-stable ring. Then $R$ is a B$\acute{e}$zout ring if and only if $R$ is
an elementary divisor ring.} \vskip2mm\hspace{-1.8em} {\it
Proof.}\ \ $\Longrightarrow $ In view of Proposition 3.1, $R$ has
stable range $2$. Thus, $R$ is Hermite, by [7, Theorem 3.4]. Suppose that
$aR+bR+cR=R$ with $a,b,c\in R$. If $a\not\in J(R)$, then there
exists a $z\in R$ such that $aR+(b+cz)R=R$. If $a\in J(R)$, then
there exist some $x_1,x_2,x_3\in R$ such that $ax_1+bx_2+cx_3=1$.
Hence, $bx_2(1-ax_1)^{-1}+cx_3(1-ax_1)^{-1}=1$. Thus,
$$\big(x_2(1-ax_1)^{-1}\big)a+bx_2(1-ax_1)^{-1}+cx_3(1-ax_1)^{-1}=1-\big(x_2(1-ax_1)^{-1}\big)a\in
U(R).$$ Therefore,
$$\begin{array}{c}
\big(x_2(1-ax_1)^{-1}\big)a\big(1-\big(x_2(1-ax_1)^{-1}\big)a\big)^{-1}+\big(bx_2(1-ax_1)^{-1}+cx_3(1-ax_1)^{-1}\big)\\
\big(1-\big(x_2(1-ax_1)^{-1}\big)a\big)^{-1}=1. \end{array}$$
Hence, $\big(x_2(1-ax_1)^{-1}\big)aR+\big(x_2(1-ax_1)^{-1}b+x_3(1-ax_1)^{-1}c\big)R=R$. In
light of Theorem 1.1, $R$ is an elementary divisor ring.

$\Longleftarrow$ This is obvious, by Theorem 1.1.\hfill$\Box$

\vskip4mm In [7, Remark 4.7], W.W. McGovern asked a question: is
there any elementary divisor domain which does not have almost
stable range 1? This was affirmatively answered by Roitman in
[12, Example 3.3]. In fact, J-stable B$\acute{e}$zout domains which do not
have almost stable range 1 provide rich such examples. For
instance, construct $R$ as in Example 2.13. Then $R$ is a J-stable
B$\acute{e}$zout domain. Thus, $R$ is an elementary divisor ring,
in terms of Theorem 4.3. In this case, $R$ does not have almost
stable range 1. Further, we derive

\vskip4mm \hspace{-1.8em} {\bf Corollary 4.4 [7, Theorem 3.7].}\ \ {\it Let
$R$ have almost stable range 1. Then $R$ is a B$\acute{e}$zout ring if
and only if $R$ is an elementary divisor ring.}
\vskip2mm\hspace{-1.8em} {\it Proof.}\ \ By virtue of Theorem 2.1,
$R$ is a J-stable ring. This completes the proof, in terms of Theorem
4.3.\hfill$\Box$

\vskip4mm \hspace{-1.8em} {\bf Corollary 4.5.}\ \ {\it Every
J-adequate ring is an elementary divisor ring.}
\vskip2mm\hspace{-1.8em} {\it Proof.}\ \ Let $R$ be a J-adequate ring.
Then $R$ is a J-stable ring by Theorem 3.2. Therefore we complete the
proof from Theorem 4.3.\hfill$\Box$

\vskip4mm In view of Example 3.3, every B$\acute{e}$zout
NJ-ring is J-adequate. By Corollary 4.5, we give an affirmative
answer to [15, Question 4] for commutative rings.

\vskip4mm \hspace{-1.8em} {\bf Corollary 4.6.}\ \ Every Hermite NJ-ring is an elementary divisor ring.

\vskip4mm Recall that a ring $R$ satisfies $(N)$ provided that for
any $a,b\in R$ and $a\not\in J(R)$, there exists $m\in R$ such
that $bR+mR=R$ and if some $n\in R$, $nR+aR\neq R$ and $nR+bR=R$
implies $nR+mR\neq R$. We extend [4, Corollary 2.6] as follows.

\vskip4mm \hspace{-1.8em} {\bf Corollary 4.7.}\ \ {\it Every
B$\acute{e}$zout ring satisfying $(N)$ is an elementary divisor
ring.} \vskip2mm\hspace{-1.8em} {\it Proof.}\ \ As in the note in
[4, page 235] and Theorem 2.1, every B$\acute{e}$zout ring
satisfying $(N)$ is a J-stable ring. Therefore the result follows
by Theorem 4.3.\hfill$\Box$

\vskip4mm We say that $a\in R$ is $\pi$-adequate to $b\in R$
provided that there exists some $n\in {\Bbb N}$ such that $a^n$ is
adequate to $b$. A Hermite ring is called a quasi adequate ring if
for any pair of nonzero elements, at least one of these elements
is $\pi$-adequate to the other. It is obvious that every
generalized adequate ring (including adequate ring) is quasi
adequate ring (cf. [13]). We now generalize
[13, Theorem 1] as follows.

\vskip4mm \hspace{-1.8em} {\bf Theorem 4.8.}\ \ {\it Every quasi
adequate ring is an elementary divisor ring.}
\vskip2mm\hspace{-1.8em} {\it Proof.}\ \ Let $R$ be a quasi
adequate ring. By Theorem 1.1 and [12, Theorem 2.5],
it suffices to consider the matrix  $A=\left(\begin{array}{cc}
a&0\\
b&c
\end{array}
\right)$ where $aR+bR+cR=R$. If $a=0$ then $bR+cR=R$. Hence $pb+qc=1$
for some $p,q\in R$. Thus, $$\left(\begin{array}{cc}
&1\\
1&
\end{array}
\right)A\left(\begin{array}{cc}
c&p\\
b&q
\end{array}
\right)\left(\begin{array}{cc}
&1\\
1&
\end{array}
\right)=\left(\begin{array}{cc}
1&\\
&0
\end{array}
\right).$$ If $c=0$ then $aR+bR=R$, and so $pa+qb=1$ for some $p,q\in R$. Thus,
$$\left(\begin{array}{cc}
p&q\\
-b&a
\end{array}
\right)A\left(\begin{array}{cc}
1&-qc\\
&1
\end{array}
\right)=\left(\begin{array}{cc}
1&\\
&ac
\end{array}
\right).$$
So we may assume that $a,c\neq 0$ . Since $R$ is
quasi adequate ring, at least one of the elements $a,c$ is
$\pi$-adequate to the other. Let $c$ be $\pi$-adequate to $a$, then there
exist some $n\in {\Bbb N}$ such that $c^n=rs$, where $rR+aR=R$ and
$s'R+aR\neq R$ for each non-invertible divisor $s'$ of $s$. We claim
that $(a+br)R+c^nrR=R$. Since $R$ is a Hermite ring, it is a
B$\acute{e}$zout ring. If $(a+br)R+c^nrR\neq R$, then
there exists a non-invertible element $h$ of $R$ such that
$(a+br)R+c^nrR=hR$. Write $hR+rR=kR$ for some $k\in R$. Write $h=kp$ for a $p\in R$, then
$a+br=kpq$ for a $q\in R$. Write $r=km$ for some $m\in R$. Then $a=k(pq-mb)$.
Since $rR+aR=kmR+k(pq-mb)R=R$, we see that $kR=R$, and so $hR+rR=R$. Thus, $hu+r^2v=1$ for some $u,v\in R$.
Clearly, $c^nr=ht$ for a $t\in R$. Then $r^2s=ht$; hence, $(1-hu)s=htv$. We infer that $s=h(tv+us)$. By hypothesis,
$hR+aR\neq R$.

As $rR+aR=R$, we can find some $r',a'\in R$ such that $rr'+aa'=1$.
Since $aR+bR+cR=R$, we have some $x,y,z\in R$ such that $ax+by+cz=1$. Hence,
$ax+b(rr'+aa')y+c(rr'+aa')z=1$, and so $a(x+ba'y+ca'z)+(br)(r'y)+(cr)r'z=1$. This implies that
$ax'+(br)y'+(c^nr)z'=1$ for some $x',y',z'\in R$. It follows that
$a(x'-y')+(a+br)y'+(c^nr)z'=1$, and so $aR+(a+br)R+(c^nr)R=R$. This implies that
$aR+hR=R$, a contradiction.

Therefore, $(a+br)R+c^nrR=R$, and so $(a+br)R+crR=R$.
Then for the matrix $A$ we have $ \left(\begin{array}{cc}
1&r\\
0&1
\end{array}
\right)A=\left(\begin{array}{cc}
a+br&cr\\
b&c
\end{array}
\right)=B$. It suffices to prove $B$ admits a diagonal
reduction.
Write $(a+br)x+cry=1$ for some $x,y\in R$. Then the matrix  $ \left(\begin{array}{cc}
x&-cr\\
y&a+br
\end{array}
\right)$ is invertible, and we see that  $$ \left(\begin{array}{cc}
1&0\\
-(bx+cy)&1
\end{array}
\right) B   \left(\begin{array}{cc}
x&-cr\\
y&a+br
\end{array}
\right)=\left(\begin{array}{cc}
1&0\\
0&ca
\end{array}
\right),$$ as desired. Now we have, $$ \left(\begin{array}{cc}
1&0\\
-(bx+cy)&1
\end{array}
\right)  \left(\begin{array}{cc}
1&r\\
0&1
\end{array}
\right)A  \left(\begin{array}{cc}
x&-cr\\
y&a+br
\end{array}
\right)=\left(\begin{array}{cc}
1&0\\
0&ca
\end{array}
\right).$$ As $\left(\begin{array}{cc}
1&r\\
0&1
\end{array}
\right)$ and $ \left(\begin{array}{cc}
1&0\\
-(bx+cy)&1
\end{array}
\right)$ are invertible, then the multiplication of them is invertible. Hence $A$ admits a diagonal reduction.

If $a$ is $\pi$-adequate to $c$, then by similar way we have $arR+(br+c)R=R$. Write $arx+(br+c)y=1$. In this case we have  $$\left(\begin{array}{cc}
x&y\\
-(br+c)&ar
\end{array}
\right)A\left(\begin{array}{cc}
r&1\\
1&0
\end{array}
\right)\left(\begin{array}{cc}
1&-(xa+yb)\\
0&1
\end{array}
\right)=\left(\begin{array}{cc}
1&0\\
0&-ac
\end{array}
\right),$$ as required.
\hfill$\Box$

\vskip4mm \hspace{-1.8em} {\bf Corollary 4.9.}\ \ {\it Every
$\pi$-adequate ring is an elementary divisor ring.}
\vskip2mm\hspace{-1.8em} {\it Proof.}\ \ Let $R$ be a
$\pi$-adequate ring. In view of Theorem 3.8, $R$ is J-stable, and
so it has stable range 2 by Proposition 4.1. It follows from [7, Theorem 3.4]
that $R$ is Hermite. Thus, $R$ is quasi adequate. Therefore we
complete the proof, in terms of Theorem 4.8.\hfill$\Box$

\vskip4mm Furthermore, we can extend [17, Theorem 6] as follows.

\vskip4mm \hspace{-1.8em} {\bf Corollary 4.10.}\ \ {\it Every
quasi adequate ring has stable range $2$.}
\vskip2mm\hspace{-1.8em} {\it Proof.}\ \ Let $R$ be a quasi
adequate ring. Then $R$ is an elementary divisor ring, by Theorem
4.8. Hence, $R$ is a Hermite ring. Therefore the proof is true, by
[7, Theorem 3.4].\hfill$\Box$

\vskip4mm A ring $R$ is strongly completeable provided that $a_1R+\cdots +a_nR=dR, a_i,d\in R, i=2,\cdots ,n, $ implies there is a
matrix over $R$ with first row $a_1,\cdots ,a_n$ and $det(A)=d$. One easily checks that a B$\acute{e}$zout ring is strongly completable if and only
if for any $d\in R, d\in
a_1R+\cdots +a_nR, a_i\in R, i=2,\cdots , n, $ implies there is a
matrix over $R$ with first row $a_1,\cdots , a_n$ and $det(A)=d$.

\vskip4mm \hspace{-1.8em} {\bf Theorem 4.11.}\ \ {\it Every J-stable ring is strongly completable.}
\vskip2mm\hspace{-1.8em} {\it Proof.}\ \ Let $R$ be a J-stable ring. Suppose that $a_1R+\cdots +a_nR=dR, a_i,d\in R, i=2,\cdots , n$.
If $n=2$, $d=a_1x_1+a_2x_2$ for some
$x_1,x_2\in R$. Then $-x_2,x_1$ works as a second row. Suppose
that the assertion holds for $k<n (n\geq 3)$. Write $d=a_1x_1+\cdots +a_nx_n, a_1=dq_1,\cdots ,a_n=dq_n$ for some
$x_1,\cdots ,x_n,q_1,\cdots ,q_n\in R$. Let $c=x_1q_1+\cdots +x_nq_n-1$. Then $dc=0$. Further, we have $q_{1}R+\cdots +q_{n-1}R+(q_{n}x_n-c)R=R$.

Case I. $q_{1}R+\cdots +q_{(n-2)}R\nsubseteq J(R)$. In view of
Proposition 2.4, we can find some $z\in R$ such that
$q_{1}R+\cdots +q_{(n-2)}R+(q_{(n-1)}+(q_{n}x_n-c)z)R=R$. Hence,
$d\in dq_{1}R+\cdots
+dq_{(n-2)}R+(dq_{(n-1)}+d(q_{n}x_n-c)z))R=a_{1}R+\cdots
+a_{(n-2)}R+(a_{(n-1)}+a_{n}x_nz)R$. We infer that
$dR=a_{1}R+\cdots
+a_{(n-2)}R+(a_{(n-1)}+a_{n}x_nz)R$. By hypothesis, there exists an
$(n-1)\times (n-1)$ matrix $D$ whose first row is $a_1,\cdots
,a_{n-2},a_{n-1}+a_nx_nz$ and $det(D)=d$. Let
$$A=
\left(
\begin{array}{cc}
D&\begin{array}{c} a_{n}\\
0\\
\vdots\\
0
\end{array}\\
\begin{array}{ccc}
0&\cdots &0
\end{array}&1
\end{array}
\right) \left(
\begin{array}{cc}
I_{n-2}&{\bf 0}\\
{\bf 0}&
\begin{array}{cc}
1&0\\
-x_nz&1
\end{array}
\end{array}
\right).$$ Then $A$ is the required matrix.

Case II. $q_{1}R+\cdots +q_{n-2}R\subseteq J(R)$. Then
$q_{n-1}$ or $q_{n}x_n-c$ is not in $J(R)$. Suppose that
$q_{n-1}\not\in J(R)$. Then $(q_{1}+q_{n-1})R+q_{2}R+\cdots
+q_{n-2}R\nsubseteq J(R)$. Clearly,
$(q_{1}+q_{n-1})R+q_{2}R+\cdots
+q_{n-2}R+q_{n-1}R+(q_{n}x_n-c)R=R$. Similarly to Case I, we have
a $z\in R$ such that $(q_{1}+q_{n-1})R+q_{2}R+\cdots
+q_{n-2}R+(q_{n-1}+q_{n}(x_nz-c))R=R.$ Thus, $d\in
d(q_{1}+q_{n-1})R+\cdots
+dq_{n-2}R+(dq_{n-1}+dq_{n}x_nz)R=(a_{1}+a_{n-1})R+\cdots
+a_{(n-2)}R+(a_{(n-1)}+a_{n}x_nz)R$. It follows that
$dR=(a_{1}+a_{n-1})R+\cdots +a_{n-2}R+(a_{n-1}+a_{n}x_nz)R.$ By
the discussion in Case I, we can find a matrix $A$ whose fist row
is $a_1+a_{n-1},a_2,\cdots ,a_{n-1},a_n$ and $det(A)=d$. Let
$$B=A
\left(
\begin{array}{ccccc}
I_{n-2}&{\bf 0}\\
\begin{array}{cc}
-1&0\\
0&0
\end{array}&I_2
\end{array}
\right).$$ Then $B$ is the desired matrix. Suppose that
$q_{n}x_n-c\not\in J(R)$. Then $\big(q_{1}+(q_{n}x_n-c)\big)R+q_{12}R+\cdots
+q_{1(n-2)}R\nsubseteq J(R)$. Similarly, we prove that there
exists a matrixbbb whose first row is $a_1,\cdots,a_n$ and its
determinant is $d$.

By induction, the theorem is proved.\hfill$\Box$

\vskip4mm \hspace{-1.8em} {\bf Corollary 4.12.} \ \ Every J-adequate ring and every $\pi$-adequate ring are strongly completable.

\vskip4mm Immediately, we see that every adequate ring is strongly completable. Thus, every regular ring is strongly completable.

\vskip4mm \hspace{-1.8em} {\bf Example 4.13.}\ \ For any $n\in {\Bbb N}, {\Bbb Z}_n$ is strongly completable. As every homomorphic image of a principal ideal ring is principal ideal and ${\Bbb Z}_n$ is a homomorphic image of ${\Bbb Z}$, ${\Bbb Z}_n$ is a principal ideal ring, and then a B$\acute{e}$zout ring. One easily checks that every nonzero prime ideal of ${\Bbb Z}_n$ is a maximal ideal, and so every prime ideal of ${\Bbb Z}_n$ is contained in a unique maximal ideal. As it is a finite ring, for every element $s$ of ${\Bbb Z}_n$, we see that $Z(s)$, i.e., the set of maximal ideals containing $a$, is finite. In light of [4, Theorem 4.3], ${\Bbb Z}_n$ is adequate, and we are through.

\vskip4mm We note that Theorem 4.11 extend [9, Theorem 2.1] as well. Also every ring having almost stable range 1 is strongly completable. As a consequence, we have

\vskip4mm \hspace{-1.8em} {\bf Corollary 4.14 [8, Theorem].}\ \ Every Dedekind domain is strongly completable.

\vskip15mm \bc{\bf REFERENCES}\ec \vskip4mm {\small \re{1} H.
Chen, {\it Rings Related Stable Range Conditions}, Series in
Algebra 11, World Scientific, Hackensack, NJ, 2011.

\re{2} L. Gillman and M. Henriksen, Some remarks about elementary
divisor rings, {\it Trans. Amer. Math. Soc.}, {\bf 82}(1956),
362--365.

\re{3} M. Henriksen, On a class of regular rings that are
elementary divisor rings, {\it{Arch. Math.}}, {\bf{24}}(1973),
133--141.

\re{4} M. Larsen; W. Lewis and T. Shores, Elementary divisor rings
and finitely presented modules, {\it Trans. Amer. Math. Soc.},
{\bf 187}(1974), 231--248.

\re{5} N. Mahdou; A. Mimouni and M.A.S.
Moutui, On pm rings, rings of finite character and h-local rings,
{\it J. Algebra Appl.}, {\bf 6}(2014), 1450018 [11 pages], DOI:
10.1142/S0219498814500182.

\re{6} W.W. McGovern, Neat rings, {\it J. Pure Appl. Algebra}, {\bf 205}(2006), 243--265.

\re{7} W.W. McGovern, B\'{e}zout rings with almost stable range
$1$, {\it J. Pure Appl. Algebra}, {\bf 212}(2008), 340--348.

\re{8} M.E. Moore, A strongly complement property of Dedekind
domain, {\it Czechoslovak Math. J.}, {\bf 25}(100)(1975),
282--283.

\re{9} M. Moore and A. Steger, Some results on completability in commutative rings, {\it Pacific J. Math.},
{\bf 37}(1971), 453--460.

\re{10} W.K. Nicholson, Rings whose elements are quasi-regular or
regular, {\it Aequations Mathematicae}, {\bf 9}(1973), 64--70.

\re{11} A.M. Rahimi, {\it Some Results on Stable Range in Commutative Rings}, Ph.D. Thesis, The University of Texas at Arlington, 1993.

\re{12} M. Roitman, The Kaplansky condition and rings of almost
stable range $1$, {\it Proc. Amer. Math. Soc.}, {\bf 141}(2013),
3013--3018.

\re{13} B.V. Zabavsky, Generalized adequate rings, {\it Ukrainian
Math. J.}, {\bf 48}(1996), 614--617.

\re{14} B.V. Zabavsky, Diagonal Reduction of Matrices over Rings,
Mathematical Studies Monograph Series, Vol. XVI, VNTL Publisher, 2012.

\re{15} B.V. Zabavsky, Questions related to the K-theoretical
aspect of Bezout rings with various stable range conditions,
{\it Math. Stud.}, {\bf 42}(2014), 89-103.

\re{16} B.V. Zabavsky and S.I. Bilavska,
Every zero adequate ring is an exchange ring, {\it J. Math. Sci.},
{\bf 187}(2012), 153--156.

\re{17} B.V. Zabavsky and M.Y. Komarnyts'kyi, Cohn-type theorem
for adequacy and elementary divisor rings, {\it J. Math. Sci.},
{\bf 167}(2010), 107--111.

\end{document}